\documentclass[twocolumn]{autart}    

\usepackage{graphicx}          
\usepackage{IEEEtrantools}
\usepackage{amsmath}
\usepackage{amssymb}
\usepackage{caption}

\usepackage{color}

\begin{document}

\begin{frontmatter}

\newtheorem{Definition}{Definition}
\newtheorem{Lemma}{Lemma}
\newtheorem{Theorem}{Theorem}
\newtheorem{Corollary}{Corollary}
\newtheorem{proposition}{Proposition}
\newtheorem{Assumption}{Assumption}


\title{Control Barrier Functions for Stochastic Systems\thanksref{footnoteinfo}} 

\thanks[footnoteinfo]{This paper was not presented at any IFAC 
meeting. Corresponding author A. Clark.}

\author[WPI]{Andrew Clark}\ead{aclark@wpi.edu}    

\address[WPI]{Dept. of Electrical and Computer Engineering, Worcester Polytechnic Institute, 100 Institute Road, Worcester, MA 01609}  


\begin{abstract}                          
 Control Barrier Functions (CBFs) aim to ensure safety by constraining the control input at each time step so that the system state remains within a desired safe region. This paper presents a framework for CBFs in stochastic systems  in the presence of Gaussian process and measurement noise. We first consider the case where the system state is known  at each time step, and present reciprocal and zero CBF constructions that guarantee safety with probability 1. We extend our results to high relative degree systems and present both general constructions and the special case of linear dynamics and affine safety constraints. We then develop CBFs for incomplete state information environments, in which the state must be estimated using sensors that are corrupted by Gaussian noise. We prove that our proposed CBF ensures safety with probability 1 when the state estimate is within a given bound of the true state, which can be achieved using an Extended Kalman Filter when the system is linear or the process and measurement noise are sufficiently small. We propose control policies that combine these CBFs with Control Lyapunov Functions in order to jointly ensure  safety and stochastic stability. Our results are validated via numerical study on a multi-agent collision avoidance scenario.
\end{abstract}

\end{frontmatter}

\section{Introduction}
\label{sec:intro}
Safety, defined as ensuring that the state of a control system remains within a particular region, is an essential property in applications including transportation, medicine, and energy. The need for safety  has motivated extensive research into synthesizing and verifying controllers to satisfy safety requirements. Methodologies include barrier methods \cite{prajna2007framework}, discrete approximations \cite{chutinan1999hybrid,ratschan2005safety,mitra2013verifying}, and reachable set computation \cite{girard2006efficient,althoff2011reachable}. 

\begin{figure}[!ht]
\centering
\includegraphics[width=2.5in]{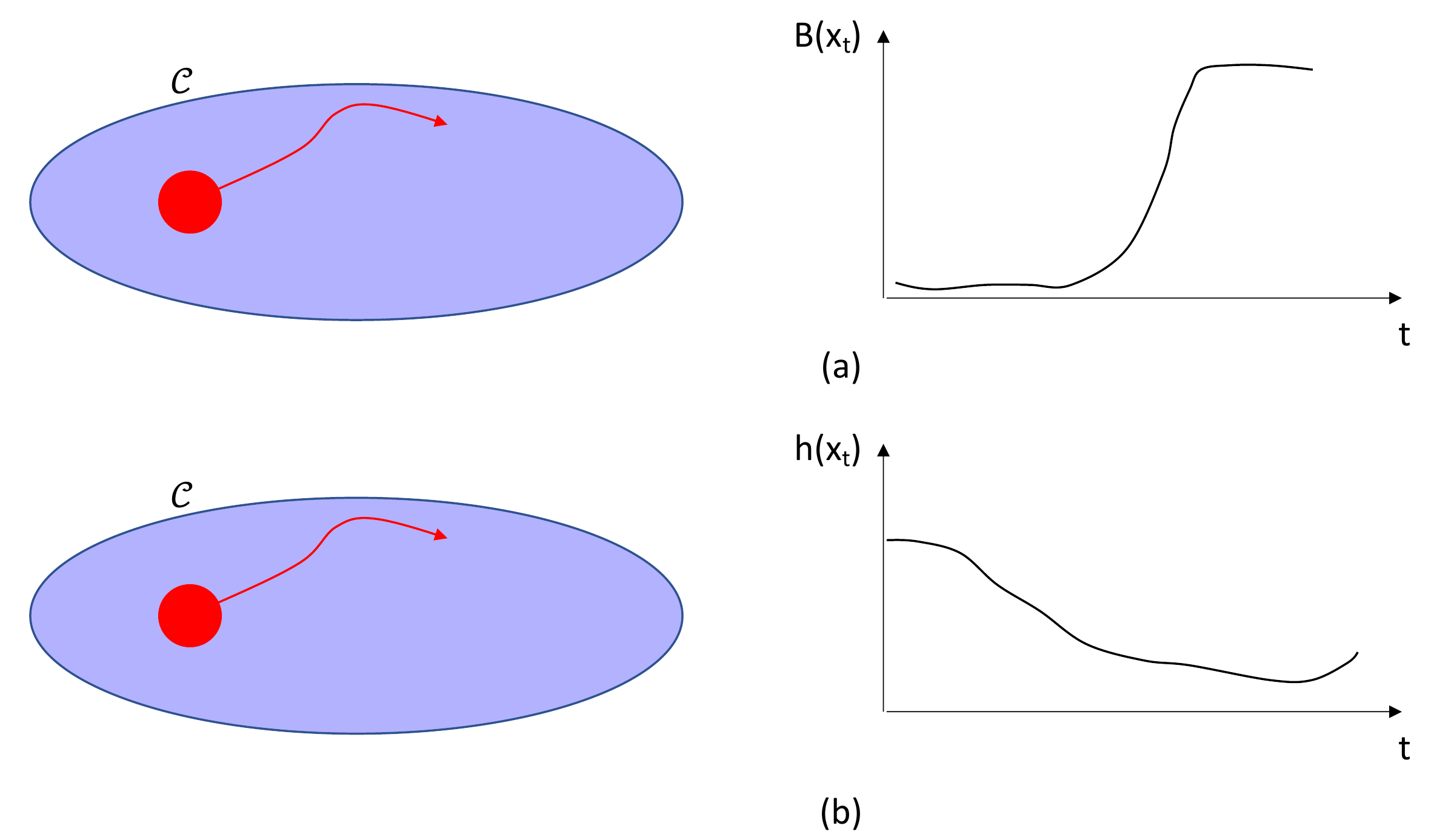}
\caption{Illustration of (a) Reciprocal CBF and (b) Zero CBF.}
\label{fig:intro}
\end{figure}

Recently, Control Barrier Functions (CBFs) have emerged as a promising approach to ensure safety while maintaining computational tractability \cite{ames2016control}. A CBF is a function that either decays to zero (Zero CBF, or ZCBF) or diverges to infinity (Reciprocal CBF, or RCBF) as the state trajectory approaches the boundary of the safe region. Safety of the system can be guaranteed by adding a constraint to the control input, which ensures that the CBF remains finite in the case of RCBF and positive in the case of ZCBF (Fig. \ref{fig:intro}). The CBF approach has been successfully applied to bipedal locomotion \cite{hsu2015control,nguyen2016optimal}, automotive control \cite{mehra2015adaptive,chen2017obstacle}, and UAVs \cite{wu2016safety}. Furthermore, by composing a CBF with a Control Lyapunov Function (CLF), optimization-based control policies with joint guarantees on safety and stability can be designed.

Existing CBF techniques are applicable to deterministic systems with exact observation of the system state. Many control systems, however, operate in the presence of noise in both the system dynamics and sensor measurements. A CBF framework for stochastic systems would enable computationally tractable control with probabilistic guarantees on safety by making the CBF method applicable to a broader class of systems.

In this paper, we generalize CBF-based methods for safe control to stochastic systems. We consider complete information systems, in which the exact system state is known, as well as incomplete information systems in which only noisy measurements of the state are available. For both cases, we formulate stochastic versions ZCBF and RCBF, and show that a linear constraint on the control at each time step results in provable safety guarantees. We make the following  contributions:

\begin{itemize}
\item In the complete information case, we formulate ZCBFs and RCBFs and derive sufficient conditions for the system to satisfy safety with probability 1.
\item In the incomplete information case, we consider a class of controllers in which the state estimate is obtained via Extended Kalman Filter (EKF). We derive bounds on the probability of violating the safety constraints as a function of the estimation error of the filter.
\item We derive sufficient conditions for constructing ZCBFs for high relative degree  systems, and analyze the special case of linear systems  with affine safety constraints and complete state information.
\item We construct optimization-based controllers that integrate stochastic CLFs with CBFs to ensure safety and performance. The controllers solve quadratic programs at each time step and thus can be implemented on embedded systems.
\item We evaluate our approach via numerical study on multi-agent collision avoidance. We find that the proposed ZCBF guarantees safety while still allowing the agents to reach their desired final states.
\end{itemize}
The rest of the paper is organized as follows. Section \ref{sec:preliminaries} presents needed background. Section \ref{sec:complete} presents CBF constructions in the complete information case. Section \ref{sec:incomplete} considers the incomplete information case. Section \ref{sec:policy} presents control policy constructions via stochastic CBFs. Section \ref{sec:simulation} contains numerical results. Section \ref{sec:conclusion} concludes the paper.
\section{Related Work}
\label{sec:related}
The CBF method for synthesizing safe controllers was proposed in \cite{ames2014control,ames2016control}. For a comprehensive survey of recent work on CBFs, see \cite{ames2019control}.  Composition of CBFs with CLFs for guaranteed safety and stability was proposed in \cite{romdlony2016stabilization}.  CBFs have been proposed for input-constrained systems \cite{rauscher2016constrained},  systems with delays \cite{jankovic2018control},  self-triggered systems~\cite{yang2019self}, and linearizable systems~\cite{xu2018constrained}. Extensions to incorporate signal temporal logic constraints were developed in \cite{lindemann2019control}. A framework for \emph{exponential CBFs} that enable safety guarantees in high relative-degree systems was proposed in \cite{nguyen2016exponential}. High relative-degree deterministic systems were also considered in \cite{xiao2019control,khojasteh2020probabilistic}. While the present paper also considers high relative degree systems, we propose a different approach and, moreover, consider the problem in a stochastic setting. Learning-based methods for CBFs in systems with incomplete information due to uncertainties were presented in \cite{cheng2019end,fan2020bayesian,khojasteh2020probabilistic,cheng2020safe}.

The problem of verifying safety of a given system and controller has been studied extensively over the past several decades \cite{chutinan1999hybrid,ratschan2005safety,dolginova1997safety,tabuada2009verification,tomlin1998conflict}.  In the verification literature, the approach that is closest to the present work is the barrier function method~\cite{prajna2004safety,prajna2007framework}. Barrier certificates provide provable guarantees that a system with given controller does not enter an unsafe region. More recently, a tighter barrier function construction that enables controller synthesis for stochastic systems was proposed in \cite{santoyo2019verification}. A discrete-time barrier certificate for ensuring satisfaction of temporal logic properties was proposed in \cite{jagtap2020formal}.  Barrier certificate methods, however, enable safety verification of a given system, but do not provide an approach for synthesizing controllers with safety guarantees. Indeed, existing techniques for synthesizing barrier certificates using sum-of-squares optimization are inapplicable to designing control barrier functions.


The preliminary conference version of this paper~\cite{clark2019control} introduced CBFs for stochastic systems, including what this paper refers to as \emph{reciprocal CBFs}. The present paper introduces the additional notion of \emph{zero CBFs} for stochastic systems, as well as methodologies for computing CBFs for high relative degree systems.  We also extend our results in the incomplete information case to systems where the output is nonlinear in the input.


\section{Background}
\label{sec:preliminaries}
This section provides background on martingales and stochastic differential equations (SDEs). In what follows, we let $\cdot^{+} = \max{\{\cdot, 0\}}$, $\cdot^{-} = \min{\{\cdot,0\}}$, $\mathbf{E}(\cdot)$ denote expectation, and $\mathbf{tr}(\cdot)$ denote the trace. A function $\alpha: \mathbb{R} \rightarrow \mathbb{R}$ is \emph{class-K} if it is strictly increasing and $\alpha(0) = 0$. We let $[x]_{i}$ denote the $i$-th element of vector $x$. 

We consider stochastic processes with respect to a probability space $(\Omega, \mathcal{F}, Pr)$, where $\Omega$ is a sample space, $\mathcal{F}$ is a $\sigma$-field over $\Omega$, and $Pr: \mathcal{F} \rightarrow [0,1]$ is a probability measure. A \emph{filtration} $\{\mathcal{F}_{t} : t \geq 0\}$ is a collection of sub-$\sigma$-fields with $\mathcal{F}_{s} \subseteq \mathcal{F}_{t} \subseteq \mathcal{F}$ for $0 \leq s < t < \infty$. A stochastic process is \emph{adapted} to filtration $\{\mathcal{F}_{t}\}$ if, for each $t \geq 0$, $X_{t}$ is an $\mathcal{F}_{t}$-measurable random variable \cite{karatzas2012brownian}.

\begin{Definition}
\label{def:martingale}
The random process $x_{t}$ is a \emph{martingale} if $\mathbf{E}(x_{t} | x_{s}) = x_{s}$ for all $t \geq s$, a \emph{submartingale} if $\mathbf{E}(x_{t} | x_{s}) \geq x_{s}$ for all $t \geq s$, and a \emph{supermartingale} if $\mathbf{E}(x_{t} | x_{s}) \leq x_{s}$ for all $t \geq s$.
\end{Definition}

A stopping time is defined as follows.

\begin{Definition}
\label{def:stopping-time}
A random time $\tau$ is a stopping time of a filtration $\mathcal{F}_{t}$ if the event $\{\tau \leq t\}$ belongs to the $\sigma$-field $\mathcal{F}_{t}$ for all $t \geq 0$.
\end{Definition}

Let $x_{t}$ be a submartingale (resp. supermartingale) and let $\tau$ be a stopping time. If $t \wedge \tau$ denotes the minimum of $t$ and $\tau$, then $x_{t \wedge \tau}$ is a submartingale (resp. supermartingale). The following result gives bounds on the maximum value of a submartingale.

\begin{Theorem}[Doob's Martingale Inequality \cite{karatzas2012brownian}]
\label{theorem:DMI}
Let $x_{t}$ be a submartingale, $[t_{0},t_{1}]$ a subinterval of $[0,\infty)$, and $\lambda > 0$. Then 
\begin{equation}
\label{eq:DMI-1}
\lambda Pr\left(\sup_{t_{0} \leq t \leq t_{1}}{x_{t}} \geq \lambda\right) \leq \mathbf{E}(x_{t_{1}}^{+}). 
\end{equation}
\end{Theorem}
The following result follows directly from Doob's Martingale Inequality.
\begin{Corollary}
\label{corollary:DMI}
Let $x_{t}$ be a supermartingale, $[t_{0},t_{1}]$ a subinterval of $[0,\infty)$, and $\lambda > 0$. Then 
\begin{equation}
\label{eq:DMI-supermart-1}
\lambda Pr\left(\inf_{t \in [t_{0},t_{1}]}{x_{t}} \leq -\lambda\right) \leq \mathbf{E}(x_{t_{1}}^{+}) - \mathbf{E}(x_{t_{1}}).
\end{equation}
\end{Corollary}
\emph{Proof:} Since $x_{t}$ is a supermartingale, $-x_{t}$ is a submartingale. Applying (\ref{eq:DMI-1}) with the submartingale $-x_{t}$ completes the proof.$\qed$ 

The \emph{quadratic variation} $\langle X \rangle$ of a random process $X$ is the unique adapted increasing process for which $\langle X \rangle_{0} = 0$ and $X^{2}-\langle X \rangle$ is a martingale \cite{karatzas2012brownian}.

We next define a semimartingale and give a composition result on semimartingales.
\begin{Definition}
\label{def:semimartingale}
A continuous semimartingale $x_{t}$ is a stochastic process which has decomposition $x_{t} = x_{0} + M_{t} + A_{t}$ with probability $1$, where $M_{t}$ is a martingale and $A_{t}$ is the difference between two continuous, nondecreasing, adapted processes. 
\end{Definition}
For any stopping time $\tau$ and semimartingale $x_{t}$, $x_{t \wedge \tau}$ is a semimartingale. The following lemma gives a composition rule for semimartinigales.
\begin{Lemma}[It\^{o}'s Lemma \cite{karatzas2012brownian}]
Let $f(x,t)$ be a twice-differentiable function and let $x_{t}$ be a semimartingale. Then $f(x_{t})$ is a semimartingale that satisfies 
\begin{multline*}
f(x_{t}) = f(x_{0}) + \int_{0}^{t}{f^{\prime}(x_{s}) \ dM_{s}} + \int_{0}^{t}{f^{\prime}(x_{s}) \ dA_{s}} \\
+ \frac{1}{2}\int_{0}^{t}{f^{\prime\prime}(x_{s}) \ d \langle M \rangle_{s}}
\end{multline*}
with probability 1 for all $t$.
\end{Lemma}

A stochastic differential equation (SDE) in It\^{o} form is defined by 
\begin{equation}
\label{eq:sde}
dx_{t} = a(x,t) \ dt + \sigma (x,t) \ dW_{t}
\end{equation}
where $a(x,t)$ and $\sigma(x,t)$ are continuous functions and $W_{t}$ is a Brownian motion. The dimension of $x_{t}$ is equal to $n$, while the dimension of $W_{t}$ is equal to $q$. A strong solution to an SDE is defined as follows.
\begin{Definition}
\label{def:strong-solution}
A strong solution of SDE (\ref{eq:sde}) with respect to Brownian motion $W_{t}$ and initial condition $\chi$ is a process $\{x_{t} : t \in [0,\infty)\}$ with continuous sample paths and the following properties:
\begin{enumerate}
\item[(i)] $Pr(x_{0} = \chi) = 1$ 
\item[(ii)] For every $1 \leq i \leq n$, $1 \leq j \leq r$, and $t \in [0,\infty)$, $$Pr\left(\int_{0}^{t}{|a_{i}(x_{\tau},\tau)| + \sigma_{ij}^{2}(x_{\tau},\tau) \ d\tau}< \infty\right) = 1.$$
\item[(iii)] The integral equation $$x_{t} = x_{0} + \int_{0}^{t}{a(\tau,x_{\tau}) \ d\tau} + \int_{0}^{t}{\sigma(\tau,x_{\tau}) \ dW_{\tau}},$$ 
where the latter term is a stochastic integral with respect to the Brownian motion $W_{t}$, holds with probability 1.
\end{enumerate}
\end{Definition}
Any strong solution of an SDE is a semimartingale. For such strong solutions, if $f(x,t)$ is a twice differentiable function  and $z_{t} = f(x_{t},t)$, then It\^{o}'s Lemma reduces to 
\begin{multline}
\label{eq:Ito-differential}
dz_{t} = \\
\left(\frac{\partial f}{\partial t} + \frac{\partial f}{\partial x}a(x,t) + \frac{1}{2}\mathbf{tr}\left(\sigma(x,t)^{T}\frac{\partial^{2} f}{\partial x^{2}}\sigma(x,t)\right)\right) \ dt \\
+ \left(\frac{\partial f}{\partial x}\sigma(x,t)\right) \ dW_{t}
\end{multline}

\section{Complete-Information CBFs}
\label{sec:complete}
This section presents our construction of control barrier functions for stochastic systems where the controller has complete state information. 

\subsection{Problem Statement}
\label{subsec:statement}
We consider a system with time-varying state $x_{t} \in \mathbb{R}^{n}$ and control input $u_{t} \in \mathbb{R}^{m}$. The state $x_{t}$ follows the SDE 
\begin{equation}
\label{eq:control-SDE}
dx_{t} = (f(x_{t}) + g(x_{t})u_{t}) \ dt + \sigma(x_{t}) \ dW_{t}
\end{equation}
where $f: \mathbb{R}^{n} \rightarrow \mathbb{R}^{n}$, $g: \mathbb{R}^{n} \rightarrow \mathbb{R}^{n \times m}$, and $\sigma : \mathbb{R}^{n} \rightarrow \mathbb{R}^{q}$ are locally Lipschitz continuous functions and $W_{t}$ is a Brownian motion. We assume that (\ref{eq:control-SDE}) has a strong solution for any control signal $u_{t}$. 

The system is required to satisfy a safety constraint for all time $t$, which is expressed as $x_{t} \in \mathcal{C}$ for all $t$ where $\mathcal{C}$ is a safe operating region. The set $\mathcal{C}$ is defined by a locally Lipschitz function $h: \mathbb{R}^{n} \rightarrow \mathbb{R}$ as 
\begin{displaymath}
\mathcal{C} = \{x : h(x) \geq 0\}, \quad \partial \mathcal{C}  = \{x : h(x) = 0 \}.
\end{displaymath}
 The set of interior points of $\mathcal{C}$ is denoted as $\mbox{int}(\mathcal{C})$.

\noindent \textbf{Problem studied:} How to design a control policy that maps the sequence $\{x_{t^{\prime}} : t^{\prime} \in [0,t)\}$ to an input $u_{t}$ such that $x_{t} \in \mathcal{C}$ for all $t$ with probability 1?

We observe that, for systems where it is not possible to design a policy that ensures safety with probability 1, there may be policies that provide policy with some probability $\epsilon \in (0,1)$. Constructing such policies is a direction for future work.

\subsection{Reciprocal Control Barrier Function Construction}
\label{subsec:RCBF}
We present our first stochastic CBF construction, which is a \emph{reciprocal} CBF (RCBF) analogous to \cite{ames2016control}. 
\begin{Definition}
\label{def:reciprocal-CBF}
Let $x_{t}$ be a stochastic process described by (\ref{eq:control-SDE}). A reciprocal CBF is a function $B: \mathbb{R}^{n} \rightarrow \mathbb{R}$ that is locally Lipschitz, twice differentiable on $\mbox{int}(\mathcal{C})$, and satisfies the following properties:
\begin{enumerate}
\item There exist class-K functions $\alpha_{1}$ and $\alpha_{2}$ such that 
\begin{equation}
\label{eq:RCBF-1}
\frac{1}{\alpha_{1}(h(x))} \leq B(x) \leq \frac{1}{\alpha_{2}(h(x))}
\end{equation}
for all $x \in \mbox{int}(\mathcal{C})$.
\item There exists a class-K function $\alpha_{3}$ such that, for all $x \in \mbox{int}(\mathcal{C})$, there exists $u \in \mathbb{R}^{m}$ such that 
\begin{multline}
\label{eq:RCBF-2}
\frac{\partial B}{\partial x}(f(x) + g(x)u) + \frac{1}{2}\mathbf{tr}\left(\sigma(x)^{T}\frac{\partial^{2}B}{\partial x^{2}}\sigma(x)\right) \\
\leq \alpha_{3}(h(x))
\end{multline}
\end{enumerate}
\end{Definition}
In the deterministic case \cite{ames2016control}, the reciprocal CBF construction ensures that  $B(x)$ tends to infinity as the system state approaches the boundary of the safe region $\mathcal{C}$. Definition \ref{def:reciprocal-CBF} extends this approach to the stochastic case by providing sufficient conditions for the system to remain bounded in expectation, and hence almost surely finite, as shown by the following theorem.

\begin{Theorem}
\label{theorem:reciprocal-CBF-safe}
Suppose that there exists an  RCBF $B$ for a controlled stochastic process $x_{t}$ described by (\ref{eq:control-SDE}), and at each time $t$, $u_{t}$ satisfies (\ref{eq:RCBF-2}). Then  $Pr(x_{t} \in \mathcal{C} \ \forall t) = 1$, provided that $x_{0} \in \mathcal{C}$.
\end{Theorem}

\emph{Proof:} {\color{blue} We will show that, for all $t$, $Pr(x_{t^{\prime}} \in \mathcal{C} \ \forall t^{\prime} < t) = 1$, and hence $$Pr(x_{t} \in \mathcal{C} \ \forall t) = \lim_{t \rightarrow \infty}{Pr(x_{t^{\prime}} \in \mathcal{C} \ \forall t^{\prime} \in [0,t])} = 1.$$} Let $B$ be a RCBF and define $B_{t} = B(x_{t})$. Since each sample path of $x_{t}$ is continuous, each sample path of $B_{t}$ is continuous. Hence, if $x_{t} \notin \mathcal{C}$ for some $t$, then there exists $t^{\prime} < t$ such that $h(x_{t^{\prime}}) = 0$ and thus $B_{t^{\prime}} = \infty$ by (\ref{eq:RCBF-1}). As a result, if for all $t > 0$ and for all $\delta \in (0,1)$, we have $$Pr\left(\sup_{t^{\prime} < t}{B_{t^{\prime}}} = \infty\right) < \delta,$$ then $Pr(x_{t} \in \mathcal{C}) = 1$ for all $t$. Equivalently, $Pr(x_{t} \in \mathcal{C}) = 1$ for all $t$ if, for all $t > 0$ and $\delta \in (0,1)$, we can construct $K > 0$ such that $Pr(\sup_{t^{\prime} < t}{B_{t^{\prime}}} = \infty) \leq Pr(\sup_{t^{\prime} < t}{B_{t^{\prime}}} > K) < \delta$. 

We construct such a $K$ as follows. 
Let $L = B_{0}$, and choose $K$ such that $$K > \frac{L + t\alpha_{3}(\alpha_{2}^{-1}(\frac{1}{L}))}{\delta}.$$ Define stopping time $\beta$ as $\beta = \inf{\{t : B_{t} = 2K\}}$. We have that $x_{t \wedge \beta}$ is a semimartingale and $x_{t \wedge \beta} \in \mbox{int}(\mathcal{C})$ for all $t$\footnote{We consider the process $x_{t \wedge \beta}$ instead of $x_{t}$ in order to ensure that $B(x_{t \wedge \beta})$ remains finite with probability 1, and hence It\^{o}'s Lemma is applicable.}. The function $B(x)$ is twice differentiable on $\mbox{int}(\mathcal{C})$, and therefore for any $x$ in a sample path of $x_{t \wedge \beta}$. Hence we can apply It\^{o}'s Lemma to obtain
\begin{IEEEeqnarray}{rCl}
\IEEEeqnarraymulticol{3}{l}{
B_{t \wedge \beta} = B_{0} + \int_{0}^{t \wedge \beta}{\left[\frac{\partial B}{\partial x}(f(x_{\tau}) + g(x_{\tau})u_{\tau})\right.}} \\
\nonumber
\label{eq:RCBF-proof-1}
&&+ \left.\frac{1}{2}\mathbf{tr}\left(\sigma(x_{\tau})^{T}\frac{\partial^{2}B}{\partial x^{2}}\sigma(x_{\tau})\right)\right] \ d\tau 
+ \int_{0}^{t \wedge \beta}{\frac{\partial B}{\partial x}\sigma(x_{\tau}) \ dW_{\tau}}
\end{IEEEeqnarray}
with probability 1. We construct a sequence of stopping times $\eta_{i}$ and $\zeta_{i}$ as 
\begin{IEEEeqnarray}{rCl}
\label{eq:up-down-1}
\eta_{0} &=& 0, \zeta_{0} = \inf{\{t : B_{t} < L\}} \\
\label{eq:up-down-2}
\eta_{i} &=& \inf{\{t : B_{t} > L, t > \zeta_{i-1}\}}, i=1,2,\ldots, \\
\label{eq:up-down-3}
\zeta_{i} &=& \inf{\{t: B_{t} < L, t > \eta_{i}\}}, i=1,2,\ldots,
\end{IEEEeqnarray}
The times $\eta_{i}$ and $\zeta_{i}$ are the up- and down-crossings of $B_{t}$ over $L$. Define a random process $\tilde{B}_{t}$ by 
\begin{multline*}
\tilde{B}_{t} = L + \sum_{i=0}^{\infty}{\left[\int_{\eta_{i} \wedge t}^{\zeta_{i} \wedge t}{\alpha_{3}(\alpha_{2}^{-1}(\frac{1}{L})) \ d\tau} \right.} \\
\left. + \int_{\eta_{i} \wedge t}^{\zeta_{i} \wedge t}{\frac{\partial B}{\partial x}\sigma(x_{\tau}) \ dW_{\tau}} \right]
\end{multline*}
We will show that, for any sample path where (\ref{eq:RCBF-proof-1}) holds, we have $B_{t \wedge \beta} \leq \tilde{B}_{t \wedge \beta}$, or equivalently, $B_{t \wedge \beta} \leq \tilde{B}_{t \wedge \beta}$ with probability 1. The proof is by induction. At time $t=0$, $B_{0} = \tilde{B}_{0} = L$. For $t \in (\eta_{i}, \zeta_{i}]$, 
\begin{IEEEeqnarray}{rCl}
\nonumber
B_{t} &=& B_{\eta_{i}} + \int_{\eta_{i}}^{t}{\left[\frac{\partial B}{\partial x}(f(x_{\tau}) + g(x_{\tau})u_{\tau})\right.} \\
\nonumber
&& \left.+ \frac{1}{2}\mathbf{tr}\left(\sigma(x_{\tau})^{T}\frac{\partial^{2}B}{\partial x^{2}}\sigma(x_{\tau})\right) \right] \ d\tau \\
\label{eq:RCBF-proof-2}
&& + \int_{\eta_{i}}^{t}{\frac{\partial B}{\partial x}\sigma(x_{\tau}) \ dW_{\tau}} \\
\label{eq:RCBF-proof-3}
\tilde{B}_{t} &=& \tilde{B}_{\eta_{i}} + \int_{\eta_{i}}^{t}{\alpha_{3}(\alpha_{2}^{-1}(\frac{1}{L})) \ d\tau} + \int_{\eta_{i}}^{t}{\frac{\partial B}{\partial x}\sigma(x_{\tau}) \ dW_{\tau}} \IEEEeqnarraynumspace
\end{IEEEeqnarray}
By induction, $B_{\eta_{i}} \leq \tilde{B}_{\eta_{i}}$. The third terms of (\ref{eq:RCBF-proof-2}) and (\ref{eq:RCBF-proof-3}) are equal. It remains to show that the second term of (\ref{eq:RCBF-proof-2}) is a lower bound on the second term of (\ref{eq:RCBF-proof-3}). By definition of $\eta_{i}$, $B_{\tau} \geq L$ for all $\tau \in [\eta_{i},t]$, or equivalently, $\frac{1}{B_{\tau}} \leq \frac{1}{L}$. By Eq. (\ref{eq:RCBF-1}), $B_{\tau} \leq \frac{1}{\alpha_{2}(h(x_{\tau}))}$, and hence $\alpha_{2}(h(x_{\tau})) \leq \frac{1}{B_{\tau}}$ and $h(x_{\tau}) \leq \alpha_{2}^{-1}(\frac{1}{B_{\tau}})$. Thus $h(x_{\tau}) \leq \alpha_{2}^{-1}(\frac{1}{B_{\tau}})$ and $\alpha_{3}(h(x_{\tau})) \leq \alpha_{3}(\alpha_{2}^{-1}(\frac{1}{B_{\tau}}))$. Combining these inequalities with (\ref{eq:RCBF-2}), we obtain 
\begin{multline*}
\frac{\partial B}{\partial x}(f(x_{\tau}) + g(x_{\tau})u_{\tau}) + \frac{1}{2}\mathbf{tr}\left(\sigma(x_{\tau})^{T}\frac{\partial^{2}B}{\partial x^{2}}\sigma(x_{\tau})\right) \\
\leq \alpha_{3}\left(\alpha_{2}^{-1}\left(\frac{1}{L}\right)\right),
\end{multline*}
 and therefore the integrand of the second term of (\ref{eq:RCBF-proof-2}) is a lower bound on the integrand of the second term of (\ref{eq:RCBF-proof-3}). In particular, $L = B_{\zeta_{i}} \leq \tilde{B}_{\zeta_{i}}$. 

For $t \in [\zeta_{i},\eta_{i+1}]$, 
\begin{IEEEeqnarray*}{rCl}
\tilde{B}_{t} &=& L + \sum_{j=0}^{i}{\left[\int_{\eta_{j}}^{\zeta_{j}}{\alpha_{3}(\alpha_{2}^{-1}(\frac{1}{L})) \ d\tau} + \int_{\eta_{j}}^{\zeta_{j}}{\frac{\partial B}{\partial x}\sigma(x_{\tau}) \ dW_{\tau}}\right]} \\
&=& \tilde{B}_{\zeta_{i}} \geq L \geq B_{t}
\end{IEEEeqnarray*}
by definition of $\eta_{i}$ and $\zeta_{i}$. Hence $B_{t} \leq \tilde{B}_{t}$ for all $t$ almost surely. As a corollary, $\tilde{B}_{t \wedge \beta} \geq B_{t \wedge \beta}$ almost surely, and we have 
\begin{IEEEeqnarray}{rCl}
\label{eq:RCBF-proof-4}	
Pr\left(\sup_{t^{\prime} \in [0,t]}{B_{t^{\prime}}} > K\right) &=& Pr\left(\sup_{t^{\prime} \in [0,t]}{B_{t^{\prime} \wedge \beta}} > K\right)
\IEEEeqnarraynumspace
 \\
\nonumber 
&\leq& Pr\left(\sup_{t^{\prime} \in [0,t]}{\tilde{B}_{t^{\prime}}} > K\right)
\end{IEEEeqnarray}
Eq. (\ref{eq:RCBF-proof-4}) holds since $B_{t} = B_{t \wedge \beta}$ when $t < \beta$,  and hence, if $B_{t^{\prime}} > K$ for some $t^{\prime} < t$, then $B_{t^{\prime} \wedge \beta} > K$. It therefore suffices to prove that $Pr(\sup_{t^{\prime} < t}{\tilde{B}_{t^{\prime} \wedge \beta}} > K) < \delta$. We first show that $\tilde{B}_{t}$ is a submartingale. We have 
\begin{IEEEeqnarray*}{rCl}
\mathbf{E}(\tilde{B}_{t} | \tilde{B}_{s}) &=& \tilde{B}_{s} + \mathbf{E}\left[\sum_{i=0}^{\infty}{\int_{\eta_{i} \wedge t}^{\zeta_{i} \wedge t}{\alpha_{3}(\alpha_{2}^{-1}(\frac{1}{L})) \ d\tau}} \right. \\
&& \left. + \int_{\eta_{i} \wedge t}^{\zeta_{i} \wedge t}{\frac{\partial B}{\partial x}\sigma(x_{\tau}) \ dW_{\tau}} \right] \\
&=& \tilde{B}_{s} + \mathbf{E}\left[\sum_{i=0}^{\infty}{\int_{\eta_{i} \wedge t}^{\zeta_{i} \wedge t}{\alpha_{3}(\alpha_{2}^{-1}(\frac{1}{L})) \ d\tau}}\right] \geq \tilde{B}_{s}
\end{IEEEeqnarray*}
implying that $\tilde{B}_{t}$ is a submartingale. 

Doob's Martingale Inequality (Theorem \ref{theorem:DMI}) then yields 
\begin{IEEEeqnarray*}{rCl}
KPr\left(\sup_{\tau \in [0,t]}{\tilde{B}_{\tau \wedge \beta}} > K\right) &\leq& \mathbf{E}(\tilde{B}_{t \wedge \beta}) \\
&\leq& L + \mathbf{E}(t \wedge \beta)\alpha_{3}(\alpha_{2}^{-1}(\frac{1}{L})) \\
&\leq& L + t\alpha_{3}(\alpha_{2}^{-1}(\frac{1}{L})).
\end{IEEEeqnarray*}
Rearranging terms and using the choice of $K$ implies that $$Pr\left(\sup_{\tau \in [0,t]}{B_{\tau \wedge \beta}} > K\right) < \delta,$$ as desired. \qed 

 Theorem \ref{theorem:reciprocal-CBF-safe} implies that, by choosing $u_{t}$ at each time $t$ to satisfy (\ref{eq:RCBF-2}), safety is guaranteed with probability 1. 

\subsection{Zero Control Barrier Function Construction}
\label{subsec:zero-CBF}
An alternative construction for CBFs is the \emph{zero-CBF} (ZCBF). The deterministic ZCBF ensures that $\frac{dh}{dt} = 0$ when $h(x) = 0$, so that the system does not enter the unsafe area. Ensuring that $\frac{dh}{dt} = 0$, however, may be inadequate in the presence of stochastic noise.  
We present a zero-CBF construction for stochastic systems that generalizes the construction in the deterministic case by using the It\^{o} derivative instead of the Lie derivative. 

\begin{Definition}
\label{def:zero-CBF}
The function $h(x)$ serves as a zero-CBF for a system described by SDE (\ref{eq:control-SDE}) if  for all $x$ satisfying $h(x) > 0$, there is a $u$ satisfying
\begin{equation}
\label{eq:ZCBF}
\frac{\partial h}{\partial x}(f(x) + g(x)u) + \frac{1}{2}\mathbf{tr}\left(\sigma^{T}\frac{\partial^{2}h}{\partial x^{2}}\sigma\right) \geq -h(x)
\end{equation}
\end{Definition}

We next state the main result on safety via zero-CBFs.

\begin{Theorem}
\label{theorem:ZCBF-safety}
If $u_{t}$ satisfies  (\ref{eq:ZCBF}) for all $t$, then $Pr(x_{t} \in \mathcal{C} \ \forall t) = 1$, provided $x_{0} \in \mathcal{C}$.
\end{Theorem}

\emph{Proof:}  We will show that, for all $t$, $Pr(x_{t^{\prime}} \in \mathcal{C} \ \forall t^{\prime} < t) = 1$, and hence $$Pr(x_{t} \in \mathcal{C} \ \forall t) = \lim_{t \rightarrow \infty}{Pr(x_{t^{\prime}} \in \mathcal{C} \ \forall t^{\prime} \in [0,t])} = 1.$$ It is sufficient to show that, for any $t > 0$, any $\epsilon > 0$, and any $\delta \in (0,1)$, $$Pr\left(\inf_{t^{\prime} < t}{h(x_{t^{\prime}})} < -\epsilon\right) < \delta.$$ Let $\theta = \min{\left\{\frac{\delta\epsilon}{2t}, h(x_{0})\right\}}.$ By It\^{o}'s Lemma, 
we have that $h(x_{t})$ is given by 
\begin{IEEEeqnarray}{rCl}
\nonumber
h(x_{t}) &=& h(x_{0}) + \int_{0}^{t}{\frac{\partial h}{\partial x}(f(x_{\tau}) + g(x_{\tau})u_{\tau})} \\ 
\nonumber
&&+ \frac{1}{2}\mathbf{tr}\left(\sigma(x_{\tau})^{T}\frac{\partial^{2}h}{\partial x^{2}}\sigma(x_{\tau})\right) \ d\tau \\
&&
\label{eq:zero-CBF-proof1}
 + \int_{0}^{t}{\sigma(x_{\tau})\frac{\partial h}{\partial x} \ dW_{\tau}}
 \IEEEeqnarraynumspace
\end{IEEEeqnarray}

 We construct a sequence of stopping times $\eta_{i}$ and $\zeta_{i}$ for $i=0,1,\ldots$ as 
 \begin{IEEEeqnarray*}{rCl}
 \eta_{0} &=& 0, \quad \zeta_{0} = 	\inf{\{t : h(x_{t}) > \theta\}} \\
 \eta_{i} &=& \inf{\{t : h(x_{t}) < \theta, t > \zeta_{i-1}\}}, i=1,2,\ldots, \\
 \zeta_{i} &=& \inf{\{t : h(x_{t}) > \theta, t > \eta_{i-1}\}}, i=1,2,\ldots,
 \end{IEEEeqnarray*}
The stopping times $\eta_{i}$ and $\zeta_{i}$ are the down- and up-crossings of $h(x_{t})$ over $\theta$, respectively. Define a random process $U_{t}$ as follows. Let $U_{0} = \theta$, and let $U_{t}$ be given by 
\begin{equation}
\label{eq:zero-CBF-proof2}
U_{t} = U_{0} + \sum_{i=0}^{\infty}{\left[\int_{\eta_{i} \wedge t}^{\zeta_{i} \wedge t}{-\theta \ d\tau} + \int_{\eta_{i} \wedge t}^{\zeta_{i} \wedge t}{\sigma\frac{\partial h}{\partial x} \ dW_{\tau}}\right]}.
\end{equation}
 We have that $U_{t}$ is a semimartingale. Furthermore, we have
 \begin{IEEEeqnarray*}{rCl}
 \mathbf{E}(U_{t}|U_{s}) &=& U_{s} + \mathbf{E}\left(\sum_{i=0}^{\infty}{\left[\int_{\eta_{i} \wedge t}^{\zeta_{i} \wedge t}{-\theta \ d\tau}\right.}\right. \\
 && \left. \left.+ \int_{\eta_{i} \wedge t}^{\zeta_{i} \wedge t}{\sigma \frac{\partial h}{\partial x} \ dW_{\tau}} \right] \right) \\
 &=& U_{s} + \mathbf{E}\left(\sum_{i=0}^{\infty}{\int_{\eta_{i} \wedge t}^{\zeta_{i} \wedge t}{-\theta \ d\tau}}\right) \leq U_{s}
 \end{IEEEeqnarray*}
and therefore $U_{t}$ is a supermartingale.

  We will first prove by induction that $h(x_{t}) \geq U_{t}$ and $U_{t} \leq \theta$. Initially, $U_{0} = \theta \leq h(x_{0})$ by construction. Suppose the result holds up to time $t \in [\eta_{i}, \zeta_{i}]$ for $i \geq 0$. Then the first term of (\ref{eq:zero-CBF-proof1}) is an upper bound on the first term of (\ref{eq:zero-CBF-proof2}) and the third terms are equal. For $t \in [\eta_{i},\zeta_{i}]$, $h(x_{t})$ and $U_{t}$ are given by
  \begin{IEEEeqnarray}{rCl}
  \nonumber
  h(x_{t}) &=& h(x_{\eta_{i}}) + \int_{\eta_{i}}^{t}{\left[\frac{\partial h}{\partial x}(f(x_{\tau}) + g(x_{\tau})u)\right.}\\
  \
  \label{eq:zero-CBF-proof3}
  && \left. + \frac{1}{2}\mathbf{tr}\left(\sigma^{T}\frac{\partial^{2}h}{\partial x^{2}}\sigma\right)\right] \ d\tau + \int_{\eta_{i}}^{t}{\sigma\frac{\partial h}{\partial x} \ dW_{\tau}} \\
  \label{eq:zero-CBF-proof4}
  U_{t} &=& U_{\eta_{i}} + \int_{\eta_{i}}^{t}{-\theta \ d\tau} + \int_{\eta_{i}}^{t}{\sigma\frac{\partial h}{\partial x} \ dW_{\tau}}.
  \end{IEEEeqnarray}
	We have that $U_{\eta_{i}} \leq h(x_{\eta_{i}}) = \theta$ by induction, and the third terms of (\ref{eq:zero-CBF-proof3}) and (\ref{eq:zero-CBF-proof4}) are equal. Since $h(x_{t}) \leq \theta$ for $t \in [\eta_{i},\zeta_{i}]$,   Eq. (\ref{eq:ZCBF}) implies 
 \begin{multline*}
 \frac{\partial h}{\partial x}(f(x)+g(x)u) + \frac{1}{2}\mathbf{tr}\left(\sigma^{T}\frac{\partial^{2}h}{\partial x^{2}}\sigma\right) \\
 \geq -h(x) \geq -\theta	
 \end{multline*}
so that the integrand of 
  the second term of (\ref{eq:zero-CBF-proof3}) is an upper bound on the integrand of the second term of (\ref{eq:zero-CBF-proof4}). Hence $h(x_{t}) \geq U_{t}$.  Furthermore, for $ t \in [\eta_{i},\zeta_{i}]$, $h(x_{t}) \leq \theta$, and thus $U_{t} \leq \theta$. 
  
  For $t \in [\zeta_{i},\eta_{i+1}]$, we have that $$U_{t} = U_{\zeta_{i}} \leq h(x_{\zeta_{i}}) = \theta \leq h(x_{t})$$ by definition of $\zeta_{i}$.
 
 Since $U_{t} \leq h(x_{t})$, we have that $$Pr\left(\inf_{t^{\prime} < t}{h(x_{t^{\prime}})} < -\epsilon\right) \leq Pr\left(\inf_{t^{\prime} < t}{U_{t^{\prime}}} < -\epsilon\right).$$ Corollary \ref{corollary:DMI} implies that   $$\epsilon Pr\left(\inf_{t^{\prime} < t}{U_{t^{\prime}}} < -\epsilon\right) \leq \mathbf{E}(U_{t}^{+})-\mathbf{E}(U_{t}).$$  The expectation of $U_{t}$ is bounded as follows. Taking expectation of both sides of (\ref{eq:zero-CBF-proof2}) yields $$\mathbf{E}(U_{t}) = U_{0} + \mathbf{E}\left[\sum_{i=0}^{\infty}{\int_{\eta_{i} \wedge t}^{\zeta_{i} \wedge t}{-\theta \ d\tau}}\right].$$ Since $\theta > 0$, the second term is bounded below by $-\theta t$, and so we have $\mathbf{E}(U_{t}) \geq \theta - \theta t$. 
  Since $U_{t} \leq \theta$, we have $\mathbf{E}(U_{t}^{+}) \leq \theta$. Combining these yields $$\mathbf{E}(\overline{U}_{t}^{+}) \leq \theta t - \theta + \theta = \theta t.$$ We therefore have 
 \begin{IEEEeqnarray*}{rCl}
 Pr\left(\inf_{t^{\prime} < t}{h(x_{t^{\prime}})} < -\epsilon\right) &\leq& Pr\left(\inf_{t^{\prime} < t}{U_{t^{\prime}}} < -\epsilon \right) \\
 &\leq& \frac{\theta t}{\epsilon} \leq \frac{\delta\epsilon}{2t}\frac{t}{\epsilon} < \delta,
 \end{IEEEeqnarray*}
 completing the proof.
    \qed

\subsection{High Relative Degree Systems}
\label{subsec:high-degree-complete}
The safety guarantees of the preceding section rely on the existence of a control input satisfying (\ref{eq:RCBF-2}) and (\ref{eq:ZCBF}) at each time $t$.  The conditions (\ref{eq:RCBF-2}) and (\ref{eq:ZCBF}), however, may fail if $\frac{\partial h}{\partial x}g(x) =0$.
In systems with high relative degree, however, it may be the case that $\frac{\partial h}{\partial x}g(x) = 0$ for some $x$, potentially preventing the system from satisfying the conditions and rendering the safety guarantees inapplicable. We propose an approach to constructing ZCBFs for such high-degree systems.   We define a set of functions $h_{i}(x)$ for $i=0,1, \ldots,$ as $h_{0}(x) = h(x)$,
\begin{equation}
\label{eq:h-i-def-0}
h_{i+1}(x) = \frac{\partial h_{i}}{\partial x}f(x) + \frac{1}{2}\mathbf{tr}\left(\sigma^{T}\left(\frac{\partial^{2}h_{i}}{\partial x^{2}}\right)\sigma\right) + h_{i}(x).
\end{equation}
This approach is similar to the high relative degree stochastic RCBF construction presented in \cite{sarkar2020high}, albeit for stochastic ZCBF. Define $\mathcal{C}_{i} = \{x : h_{i}(x) \geq 0\}$ and $$\overline{\mathcal{C}}_{r} = \bigcap_{i=0}^{r}{\mathcal{C}_{i}}.$$ 

\begin{Theorem}
\label{theorem:stochastic-ZCBF-high-degree}
Suppose that there exists $r$ such that, for any $x \in \overline{\mathcal{C}}_{r}$, we have $\frac{\partial h_{i}}{\partial x}g(x) \geq 0$ for $i < r$ and 
\begin{equation}
\label{eq:high-degree-general-condition}
\frac{\partial h_{r}}{\partial x}(f(x) + g(x)u) + \frac{1}{2}\mathbf{tr}\left(\sigma^{T}\frac{\partial^{2}h_{r}}{\partial x^{2}}\sigma\right) \geq - h_{r}(x).
\end{equation}
 Then $Pr(x_{t} \in \mathcal{C} \ \forall t) = 1$ if $x_{0} \in \overline{\mathcal{C}}_{r}$. 
\end{Theorem}

\emph{Proof:} Suppose that $u_{t}$ satisfying the conditions of the theorem is chosen at each time $t$. By Theorem \ref{theorem:ZCBF-safety}, (\ref{eq:high-degree-general-condition}) implies that $h_{r}(x_{t}) \geq 0$ for all $t$. By definition of $h_{r}(x)$ and the assumption that $\frac{\partial h_{r-1}}{\partial x}g(x)u \geq 0$, we also have $h_{r-1}(x_{t}) \geq 0$ for all $t$. Proceeding inductively, we then have $h_{i}(x_{t}) \geq 0$ for all $i=0,\ldots,r$, and hence in particular $h(x_{t}) = h_{0}(x_{t}) \geq 0$ for all $t$. $\qed$


In what follows, we show that the conditions of Theorem \ref{theorem:stochastic-ZCBF-high-degree} can be satisfied for an important subclass of systems, namely, controllable linear systems in which the safety constraint can be expressed as a half-plane. For such systems we have $f(x) = Fx$ and $g(x) = G$ for some matrices $F$ and $G$, and  the function $h(x) = a^{T}x-b$ for some $a \in \mathbb{R}^{n}$ and $b \in \mathbb{R}$. Since the system is controllable, we have $a^{T}F^{i}G \neq 0$ for some $i$. 
 The following lemma describes the structure of the $h_{i}$'s. 
\begin{Lemma}
\label{lemma:high-degree-2}
The function $h_{i}(x)$ can be written in the form
\begin{IEEEeqnarray*}{rCl}
h_{i}(x) &=& \sum_{r_{0},\ldots,r_{i-1}}{\beta_{i}^{r_{0},\ldots,r_{i-1},0}(a^{T}F^{0}x)^{r_{0}} \cdots (a^{T}F^{i-1}x)^{r_{i-1}}} \\
&& + a^{T}F^{i}x
\end{IEEEeqnarray*}
for some values of the coefficients $\beta_{i}^{r_{0},\ldots,r_{i}}$.
\end{Lemma}

\emph{Proof:} The proof is by induction on $i$. When $i=0$, the function can be written in the form $h_{i}(x) = a^{T}F^{0}x - b$, i.e., $\beta_{0}^{0} = -b$ and all other values of $\beta_{0}^{r_{0}}$ are zero. Inducting on $i$, we can then write 
\begin{IEEEeqnarray*}{rCl}
h_{i}(x) &=& z_{i}(x) + a^{T}F^{i}x \\
\frac{\partial h_{i}}{\partial x} &=& \sum_{j=0}^{i-1}{\theta_{ij}(x)a^{T}F^{j}} +a^{T}F^{i} \\
\frac{\partial^{2}h_{i}}{\partial x^{2}} &=& \sum_{j=0}^{i}{\sum_{l=0}^{i}{\zeta_{ijl}(x)(F^{j})^{T}aa^{T}F^{i}}}
\end{IEEEeqnarray*}
where the functions $z_{i}(x)$, $\theta_{ij}(x)$, and $\zeta_{ijl}(x)$ are polynomial in $(a^{T}F^{j}(x))$ for $j=0,\ldots(i,-1)$. We therefore have 
\begin{IEEEeqnarray*}{rCl}
\IEEEeqnarraymulticol{3}{l}{
h_{i+1}(x) = \left(\sum_{j=0}^{i-1}{\theta_{ij}(x)a^{T}F^{j}} + a^{T}F^{i}\right)Fx} \\
&& + \frac{1}{2}\mathbf{tr}\left(\sigma^{T}\left(\sum_{j=0}^{i}{\sum_{l=0}^{i}{\zeta_{ijl}(x)(F^{j})^{T}aa^{T}F^{l}}}\right)\sigma\right) + h_{i}(x) \\
&=& \sum_{j=0}^{i-1}{\theta_{ij}(x)a^{T}F^{j+1}x} + a^{T}F^{i+1}x \\
&& + \frac{1}{2}\sum_{j=0}^{i}{\sum_{l=0}^{i}{\zeta_{ijl}(x)\mathbf{tr}(\sigma^{T}(F^{j})^{T}aa^{T}F^{l}\sigma)}} + h_{i}(x)
\end{IEEEeqnarray*}
Hence $h_{i+1}(x)$ is a polynomial in $(a^{T}F^{0}x),\ldots,(a^{T}F^{i+1}x)$. Furthermore,  all terms except $a^{T}F^{i+1}x$ do not contain any powers of $(a^{T}F^{i+1}x)$, completing the proof. \qed

Define $r^{\prime} = \min{\{l: a^{T}F^{l}G \neq 0\}}$. By the preceding lemma, we have, for any $x \in \overline{\mathcal{C}} \triangleq \bigcap_{i=0}^{r^{\prime}}{\mathcal{C}_{i}}$, 
\begin{equation}
\label{eq:high-degree-linear}
\frac{\partial h_{i}}{\partial x}Gu = \left\{
\begin{array}{ll}
0, & i < r^{\prime} \\
a^{T}F^{r}Gu, & i=r^{\prime}
\end{array}
\right.
\end{equation}
We are now ready to state the safety result for high relative degree LTI systems.
\begin{Theorem}
\label{theorem:rel-degree-safety}
Let $r=r^{\prime}$. If $x_{0} \in \overline{\mathcal{C}}$ and 
\begin{equation}
\label{eq:ZCBF-high-degree}
\frac{\partial h_{r}}{\partial x}g(x)u \geq -\frac{\partial h_{r}}{\partial x}f(x)-\frac{1}{2}\mathbf{tr}\left(\sigma^{T}\frac{\partial^{2}h_{r}}{\partial x^{2}}\sigma\right)-h_{r}(x)	
\end{equation} 
for all $t$, then $Pr\left(x_{t} \in \overline{\mathcal{C}}\right) = 1$. In particular, $x_{t}$ satisfies the safety constraint $\{h(x_{t}) > 0\}$ with probability 1. 
\end{Theorem}
\emph{Proof:} By Theorem \ref{theorem:stochastic-ZCBF-high-degree}, it suffices to show that (\ref{eq:high-degree-general-condition}) holds and $\frac{\partial h_{i}}{\partial x}g(x)u_{t} \geq 0$ for all $i < r$. Eq. (\ref{eq:high-degree-general-condition}) holds by Eq. (\ref{eq:ZCBF-high-degree}), and $\frac{\partial h_{i}}{\partial x}g(x) = 0$ for $i < r$ by (\ref{eq:high-degree-linear}). \qed

\section{Incomplete Information CBFs}
\label{sec:incomplete}

This section presents CBF techniques for ensuring safety of stochastic systems with incomplete information due to noisy measurements. We  give the problem statement followed by RCBF and ZCBF constructions.

\subsection{Problem Statement}
\label{subsec:incomplete-statement}

We consider a system with time-varying state $x_{t} \in \mathbb{R}^{n}$, a control input $u_{t} \in \mathbb{R}^{m}$, and output $y_{t} \in \mathbb{R}^{p}$ described by  the SDEs
\begin{IEEEeqnarray}{rCl}
\label{eq:incomp-sde-1}
dx_{t} &=& (f(x_{t}) + g(x_{t})u_{t}) \ dt + \sigma_{t} \ dV_{t} \\
\label{eq:incomp-sde-2}
dy_{t} &=& b(x_{t},u_{t}) \ dt + \nu_{t} \ dW_{t}
\end{IEEEeqnarray}
where $V_{t}$ and $W_{t}$ are Brownian motions and $f: \mathbb{R}^{n} \rightarrow \mathbb{R}^{n}$, $g: \mathbb{R}^{n} \rightarrow \mathbb{R}^{m}$, and $b: \mathbb{R}^{n} \times \mathbb{R}^{m} \rightarrow \mathbb{R}^{p}$ are locally Lipschitz continuous functions. Define $\overline{f}(x,u) = f(x) + g(x)u$. Note that, unlike in the complete information case, we assume that $\sigma_{t}$ and $\nu_{t}$ do not depend on  $x_{t}$. 

 In the incomplete case, our CBF approaches are in two parts. First, we compute an estimate of the system state and construct a safe region for the estimated state based on the accuracy of the estimator. Second, we show that the problem reduces to a complete-information stochastic SDE on the estimated state value and apply the approaches developed in Section \ref{sec:complete}.

\begin{Definition}[\cite{reif2000stochastic}]
\label{def:uniform-detect}
The pair $\left[\frac{\partial \overline{f}}{\partial x}(x,u), \frac{\partial b}{\partial x}(x)\right]$ is uniformly detectable if there exists a bounded, matrix-valued function $\Lambda(x)$ and a real number $\rho > 0$ such that $$w^{T}\left(\frac{\partial \overline{f}}{\partial x}(x,u) + \Lambda(x) \frac{\partial b}{\partial x}(x)\right)w \leq -\rho||w||^{2}$$ for all $w$, $z$, and $x$.
\end{Definition}
Uniform detectability is a standard requirement for bounding the error of estimators such as the Extended Kalman Filter \cite{reif2000stochastic,li2016new,jazwinski2007stochastic}. Note that uniform detectability and detectability are equivalent for LTI systems.

The safety condition is defined as in Section \ref{subsec:statement}. In the incomplete information case, the problem studied is stated as, \emph{For given $\epsilon \in (0,1)$, how to design a control policy that maps the sequence $\{y_{t^{\prime}} : t^{\prime} \in [0,t)\}$ to an input $u_{t}$ at each time $t$ such that $Pr(x_{t} \in \mathcal{C} \ \forall t) \geq (1-\epsilon)$?} In other words, how to ensure that the system remains safe with a given probability $(1-\epsilon)$?

We use the Extended Kalman Filter (EKF) \cite{reif2000stochastic} as a state estimator. Let $\hat{x}_{t}$ denote the estimated value of $x_{t}$, and define matrix $A_{t}$ by $$A_{t} = \frac{\partial \overline{f}}{\partial x}(\hat{x}_{t},u_{t}).$$ Let $c_{t} = \frac{\partial b}{\partial x}(\hat{x}_{t})$, $R_{t} = \nu_{t}\nu_{t}^{T}$, $Q_{t}=\sigma_{t}\sigma_{t}^{T}$, and $P_{t}$ be equal to the solution to the Riccati differential equation $$\frac{dP}{dt} = A_{t}P_{t} + P_{t}A_{t}^{T} + Q_{t} - P_{t}c_{t}^{T}R_{t}^{-1}c_{t}P_{t}.$$ The EKF estimator is defined by the SDE 
\begin{equation}
\label{eq:ekf-sde}
d\hat{x}_{t} = f(\hat{x}_{t},u_{t}) \ dt + K_{t}(dy_{t} - c_{t}\hat{x}_{t} \ dt).
\end{equation}
where $K_{t} = P_{t}c_{t}R_{t}^{-1}$ is the Kalman filter gain. Under this approach, the estimation error $\zeta_{t} = x_{t}-\hat{x}_{t}$ evolves according to the SDE $$d\zeta_{t} = ((A_{t}-K_{t}c_{t})\zeta_{t} + n_{t}) \ dt + \Gamma_{t}\left(
\begin{array}{c}
dV_{t} \\
dW_{t}
\end{array}
\right),$$ where
\begin{eqnarray}
\label{eq:EKF-notation-1}
n_{t} &=& \phi(x_{t},\hat{x}_{t},u_{t}) - K_{t}\chi(x_{t},\hat{x}_{t}) \\
\label{eq:EKF-notation-2}
\Gamma_{t} &=& \sigma_{t}-K_{t}\nu_{t}
\end{eqnarray}
We make the following additional assumptions on the system dynamics to ensure stability of the EKF. 
\begin{Assumption}
\label{assump:incomplete}
The SDEs (\ref{eq:incomp-sde-1}) and (\ref{eq:incomp-sde-2}) satisfy:
\begin{enumerate}
\item There exist constants $q,r \in \mathbb{R}_{\geq 0}$ such that $\sigma_{t}\sigma_{t}^{T} \geq qI$ and $\nu_{t}\nu_{t}^{T} \geq rI$ for all $x$ and $t$\footnote{here ``$\leq$'' refers to inequality in the semidefinite cone.}.
\item The pair $\left[\frac{\partial \overline{f}}{\partial x}(x,u), \frac{\partial b}{\partial x}\right]$ is uniformly detectable.
\item There exist real numbers $\epsilon_{\phi}$, $k_{\phi}$, $\epsilon_{\chi}$, $k_{\chi}$ such that the functions $\phi$ and $\chi$ in (\ref{eq:EKF-notation-1}) and (\ref{eq:EKF-notation-2}) are bounded by 
\begin{eqnarray*}
||\phi(x,\hat{x},u)|| &\leq& k_{\phi}||x-\hat{x}||^{2} \\
||\chi(x,\hat{x})|| &\leq& k_{\chi}||x-\hat{x}||^{2}
\end{eqnarray*}
for  $x$, $\hat{x}$ satisfying  $||x-\hat{x}||_{2} \leq \epsilon_\phi$, $||x-\hat{x}|| \leq \epsilon_{\chi}$.
\end{enumerate}
\end{Assumption} 
The first assumption states that $\mathbf{E}(\sigma_{t}\sigma_{t}^{T}) \geq q I$ and $\mathbf{E}(\nu_{t}\nu_{t}^{T}) \geq rI$ for some $q,r$, and implies that the noise matrices are uniformly bounded below. The uniform detectability assumption ensures that all of the system modes can be observed, and that the covariance of the filter can be bounded, which is necessary for deriving error bounds. The third assumption states that the linearized approximation to $\overline{f}$ is approximately accurate in a neighborhood of $x$ and $\hat{x}$. We further assume that the initial state $x_{0}$ is known.
The following result describes the stability and accuracy of the EKF.
\begin{proposition}
\label{prop:EKF-accuracy}
Suppose that the conditions of Assumption 1 hold, and that there exists $\overline{c}$ such that $||c_{t}||_{2} \leq \overline{c}$ for all $t$. There exists $\delta > 0$ such that  if $\sigma_{t}\sigma_{t}^{T} \leq \delta I$ and $\nu_{t}\nu_{t}^{T} \leq \delta I$, then  for any $\epsilon > 0$, there exists $\gamma > 0$ with 
\begin{equation}
\label{eq:EKF-accuracy}
Pr\left(\sup_{t \geq 0}{||x_{t}-\hat{x}_{t}||_{2} }\leq \gamma\right) \geq 1-\epsilon.
\end{equation}
\end{proposition}
We make two remarks on Proposition \ref{prop:EKF-accuracy}. 
First, the accuracy guarantees of the EKF do not depend on the magnitude of the control input $u_{t}$. Second, if the system is highly nonlinear, then the constant $\delta > 0$ may be small \cite{reif2000stochastic}, rendering the results inapplicable. The following lemma considers the special case of LTI systems.

\begin{Lemma}
\label{lemma:LTI-error-bound}
Suppose that $f(x_{t}) = Fx_{t}$ and $g(x_{t}) = G$ for some matrices $F$ and $G$ such that $(F,G)$ is detectable. Let $\lambda^{\ast} = \sup_{t}{\lambda_{max}(P_{t})}$, where $\lambda_{max}(\cdot)$ denotes the maximum eigenvalue of a matrix. Let $\gamma = \sqrt{\frac{n\lambda^{\ast}}{\epsilon}}$. Then $Pr(x_{t} \in \mathcal{C} \ \forall t) \geq (1-\epsilon)$.
\end{Lemma}
The proof of this lemma appears in the preliminary conference version of this paper \cite{clark2019control} and is omitted due to space constraints.


 Define $\overline{h}(x) = \inf{\{h(\hat{x}) : ||x-\hat{x}||_{2} \leq \gamma\}}.$ We have that, if $\overline{h}(\hat{x}_{t}) \geq 0$ and $||x_{t}-\hat{x}_{t}||_{2} \leq \gamma$ for all $t$, then $h(x_{t}) \geq 0$ for all $t$. When $\overline{h}(x)$ is difficult to compute or non-differentiable, define $$\overline{h}_{\gamma} = \sup{\{h(x) : ||x-x^{0}||_{2} \leq \gamma \mbox{ for some } x^{0} \in h^{-1}(\{0\})\}}.$$ The following lemma gives a sufficient condition for safety of the incomplete information system.
\begin{Lemma}
\label{lemma:incomplete-CBF-helper}
If $||x_{t}-\hat{x}_{t}||_{2} \leq \gamma$ for all $t$ and $h(\hat{x}_{t}) > \overline{h}_{\gamma}$ for all $t$, then $x_{t} \in \mathcal{C}$ for all $t$.
\end{Lemma}

\emph{Proof:} Suppose that $x_{t} \notin \mathcal{C}$ for some $t$. Since each sample path of $x_{t}$ is continuous, we must have $h(x_{\tau}) = 0$ for some $\tau \in [0,t]$. By assumption, $||\hat{x}_{\tau} - x_{\tau}||_{2} \leq \gamma$, i.e., $\hat{x}_{\tau} \in B(x_{\tau},\gamma)$. Since $x_{\tau} \in h^{-1}(\{0\})$, we have 
\begin{IEEEeqnarray*}{rCl}
h(\hat{x}_{\tau}) &\leq& \sup{\{h(x) : ||x-x_{\tau}||_{2} \leq \gamma\}} \\
&\leq& \sup{\left\{h(x) : ||x-x^{0}||_{2} \leq \gamma \right.} \\
&& \left.\mbox{ for some } x^{0} \in h^{-1}(\{0\})\right\} \\
&=& \overline{h}_{\gamma}
\end{IEEEeqnarray*}
This contradicts the assumption that $h(\hat{x}_{\tau}) > \overline{h}_{\gamma}$ and hence we must have $x_{t} \in \mathcal{C}$ for all $t$. \qed

Combining Proposition \ref{prop:EKF-accuracy} and Lemma \ref{lemma:incomplete-CBF-helper}, we have that it suffices to select $\gamma$ such that $||x_{t}-\hat{x}_{t}||_{2}$ is bounded by $\gamma$ with probability $(1-\epsilon)$, and then design a control law such that $h(\hat{x}_{t}) > \overline{h}_{\gamma}$ for all $t$. Define $\hat{h}(x) = h(x)-\overline{h}_{\gamma}$. 

\subsection{Reciprocal CBF Approach}
\label{subsec:RCBF-incomplete}
The RCBF for incomplete information systems is described as follows.

\begin{Theorem}
\label{theorem:RCBF-incomplete}
Suppose that the conditions of Proposition \ref{prop:EKF-accuracy} are satisfied and there exists a function $B: \mathbb{R}^{n} \rightarrow \mathbb{R}$ and class-K functions $\alpha_{1}$, $\alpha_{2}$, and $\alpha_{3}$ such that 
\begin{IEEEeqnarray}{rCl}
\label{eq:RCBF-incomplete-1}
\frac{1}{\alpha_{1}(\hat{h}(x))} \leq B(x) &\leq& \frac{1}{\alpha_{2}(\hat{h}(x))} \\
\IEEEeqnarraymulticol{3}{l}{
\nonumber
\frac{\partial B}{\partial x}\overline{f}(\hat{x}_{t},u_{t}) + \gamma||\frac{\partial B}{\partial x}K_{t}c||_{2} + \frac{1}{2}\mathbf{tr}\left(\nu_{t}^{T}K_{t}^{T}\frac{\partial^{2}B}{\partial x^{2}}K_{t}\nu_{t}\right)} \\
\label{eq:RCBF-incomplete-2}
&\leq& \alpha_{3}(\hat{h}(\hat{x}_{t}))
\end{IEEEeqnarray}
and $\gamma$ satisfies (\ref{eq:EKF-accuracy}) for some $\epsilon > 0$.  $Pr(x_{t} \in \mathcal{C} \ \forall t) \geq (1-\epsilon)$ if $\hat{h}(x_{0}) > 0$.
\end{Theorem}

\emph{Proof:} We show that $\hat{h}(\hat{x}_{t}) \geq 0$ for all $t$ if $||x_{t}-\hat{x}_{t}||_{2} \leq \gamma$ for all $t$. Combining Eqs. (\ref{eq:incomp-sde-2}) and (\ref{eq:ekf-sde}) yields 
\begin{IEEEeqnarray*}{rCl}
d\hat{x}_{t} &=& \overline{f}(\hat{x}_{t},u_{t}) \ dt + K_{t}(cx_{t} \ dt + \nu_{t} \ dW_{t} - c\hat{x}_{t} \ dt) \\
&=& (\overline{f}(\hat{x}_{t},u_{t}) + K_{t}c(x_{t}-\hat{x}_{t})) \ dt + K_{t}\nu_{t} \ dW_{t}
\end{IEEEeqnarray*}
Define $B_{t} = B(\hat{x}_{t})$. Hence 
\begin{multline}
dB_{t} = \left(\frac{\partial B}{\partial x}(\overline{f}(\hat{x}_{t},u_{t}) + K_{t}c(x_{t}-\hat{x}_{t})) \right. \\
\left. + \frac{1}{2}\mathbf{tr}\left(\nu_{t}^{T}K_{t}^{T}\frac{\partial^{2}B}{\partial x^{2}}K_{t}\nu_{t}\right)\right) \ dt + \frac{\partial B}{\partial x}K_{t}\nu_{t} \ dW_{t}
\end{multline}
If $||x_{t}-\hat{x}_{t}||_{2} \leq \gamma$, then $$\frac{\partial B}{\partial x}K_{t}c(x_{t}-\hat{x}_{t}) \leq ||\frac{\partial B}{\partial x}K_{t}c||_{2}||x_{t}-\hat{x}_{t}||_{2} \leq \gamma||\frac{\partial B}{\partial x}K_{t}c||_{2}.$$
Hence, if (\ref{eq:RCBF-incomplete-2}) holds, then 
\begin{IEEEeqnarray*}{rCl}
\IEEEeqnarraymulticol{3}{l}{
\frac{\partial B}{\partial x}(\overline{f}(\hat{x}_{t},u_{t}) + K_{t}c(x_{t}-\hat{x}_{t})) + \frac{1}{2}\mathbf{tr}\left(\nu_{t}^{T}K_{t}^{T}\frac{\partial^{2}B}{\partial x^{2}}K_{t}\nu_{t}\right)} \\
&\leq& \frac{\partial B}{\partial x}\left(\overline{f}(\hat{x}_{t},u_{t}) + \gamma||\frac{\partial B}{\partial x}K_{t}c||_{2}\right) \\
&& + \frac{1}{2}\mathbf{tr}\left(\nu_{t}^{T}K_{t}^{T}\frac{\partial^{2}B}{\partial x^{2}}K_{t}\nu_{t}\right) \leq \alpha_{3}(\hat{h}(\hat{x}_{t}))
\end{IEEEeqnarray*}
and thus $Pr(\hat{h}(\hat{x}_{t}) \geq 0 \ \forall t) = 1$ by Theorem \ref{theorem:reciprocal-CBF-safe}. Hence, by Lemma \ref{lemma:incomplete-CBF-helper},  $Pr(h(x_{t}) \geq 0 \ \forall t | ||x_{t}-\hat{x}_{t}||_{2} \leq \gamma \ \forall t) = 1$, and so $Pr(h(x_{t}) \geq 0) \geq 1-\epsilon$. \qed

Theorem \ref{theorem:RCBF-incomplete} implies that, if the parameter $\gamma$ is chosen such that the estimation error remains bounded by $\gamma$ with sufficient probability, then selecting a control input $u_{t}$ at each time $t$ such that (\ref{eq:RCBF-incomplete-2}) holds is sufficient to ensure safety. This constraint is linear in $u_{t}$, and all other parameters can be evaluated based on the noise characteristics and system and Kalman filter matrices.

\subsection{Zero CBF Construction}
\label{subsec:incomplete-ZCBF}
The following definition describes the zero CBF in the incomplete information case.
\begin{Definition}
\label{def:incomplete-ZCBF}
The function $\hat{h}(x)$ serves as a zero CBF for an incomplete-information system described by (\ref{eq:incomp-sde-1}) and (\ref{eq:incomp-sde-2}) if for all $x$ satisfying $\hat{h}(x) > 0$, there exists $u$ satisfying 
\begin{multline}
\frac{\partial \hat{h}}{\partial x}g(x)u 
\label{eq:incomplete-ZCBF}
\geq -\frac{\partial \hat{h}}{\partial x}f(\hat{x}_{t}) + ||\frac{\partial \hat{h}}{\partial x}K_{t}c||_{2}\gamma \\
- \frac{1}{2}\mathbf{tr}\left(\sigma^{T}\frac{\partial^{2}\hat{h}}{\partial x^{2}}\sigma\right) - \hat{h}(\hat{x}_{t})
\end{multline}
\end{Definition}
The following theorem describes the safety guarantees of the incomplete-information ZCBF.
\begin{table*}
\centering
\caption{Constraints for the CBF-based control policies in complete and information systems (Eqs. (\ref{eq:comp-info-opt}) and (\ref{eq:incomplete-info-opt}))}
\begin{tabular}{|c|c|}
\hline 
\textbf{Solution Approach} & \textbf{Linear Constraint $\Omega_{t}$} \\
\hline
\textbf{RCBF, complete information} & $\frac{\partial B}{\partial x}g(x_{t})u_{t} \leq \alpha_{3}(h(x_{t})) - \frac{\partial B}{\partial x}f(x_{t}) - \frac{1}{2}\mathbf{tr}\left(\sigma(x_{t})^{T}\frac{\partial^{2}B}{\partial x^{2}}\sigma(x_{t})\right)$ \\
\hline
\textbf{ZCBF, complete information} & $\frac{\partial h}{\partial x}g(x_{t})u_{t} \geq -\frac{\partial h}{\partial x}f(x_{t}) - \frac{1}{2}\mathbf{tr}(\sigma^{T}\frac{\partial^{2}h}{\partial x^{2}}\sigma) - h(x_{t})$ 
\\
\hline
\textbf{RCBF, incomplete information} & $\frac{\partial B}{\partial x}g(\hat{x}_{t})u \leq \alpha_{3}(\hat{h}(\hat{x}_{t})) - \frac{\partial B}{\partial x}f(\hat{x}_{t}) - \gamma||\frac{\partial B}{\partial x}K_{t}c||_{2} - \frac{1}{2}\mathbf{tr}\left(\nu_{t}^{T}K_{t}^{T}\frac{\partial^{2}B}{\partial x^{2}}K_{t}\nu_{t}\right)$ \\
\hline
\textbf{ZCBF, incomplete information} & $\frac{\partial \hat{h}}{\partial x}g(x)u 
\geq -\frac{\partial \hat{h}}{\partial x}f(\hat{x}_{t}) + ||\frac{\partial \hat{h}}{\partial x}K_{t}c||_{2}\gamma 
- \frac{1}{2}\mathbf{tr}\left(\sigma^{T}\frac{\partial^{2}\hat{h}}{\partial x^{2}}\sigma\right) - \hat{h}(\hat{x}_{t})$  \\
\hline
\end{tabular}
\label{table:linear-constraints}
\end{table*}

\begin{Theorem}
\label{theorem:incomplete-ZCBF-safety}
Suppose that $x_{0}$ satisfies $\hat{h}(x_{0}) > 0$ and, at each time $t$, $u_{t}$ satisfies (\ref{eq:incomplete-ZCBF}). If the conditions of Proposition \ref{prop:EKF-accuracy} are satisfied, then  $Pr\left(x_{t} \in \mathcal{C} \ \forall t\right) \geq (1-\epsilon)$.
\end{Theorem}
\emph{Proof:} Our approach is to show that $\hat{h}(\hat{x}_{t}) > 0$ for all $t$ when $||x_{t}-\hat{x}_{t}||_{2} \leq \gamma$, and hence safety is satisfied with probability at least $(1-\epsilon)$ by Lemma \ref{lemma:incomplete-CBF-helper}. The dynamics of $\hat{x}_{t}$ are given by the SDE (\ref{eq:incomp-sde-1}). Note that 
\begin{equation}
\label{eq:incomplete-ZCBF-1}
-\frac{\partial \hat{h}}{\partial \hat{x}}K_{t}c(x_{t}-\hat{x}_{t}) \leq ||\frac{\partial \hat{h}}{\partial \hat{x}}K_{t}c||_{2}||x_{t}-\hat{x}_{t}||_{2} \leq ||\frac{\partial \hat{h}}{\partial \hat{x}}K_{t}c||_{2}\gamma.
\end{equation}
We then have 
\begin{IEEEeqnarray}{rCl}
\IEEEeqnarraymulticol{3}{l}{
-\frac{\partial \hat{h}}{\partial \hat{x}}\left(f(\hat{x}_{t}) + K_{t}c(x_{t}-\hat{x}_{t})\right) 
 -\frac{1}{2}\mathbf{tr}\left(\sigma^{T}\left(\frac{\partial^{2}\hat{h}}{\partial \hat{x}^{2}}\right)\sigma\right)}
 \IEEEeqnarraynumspace
  \\
\label{eq:incomplete-ZCBF-2}
&\leq& -\frac{\partial \hat{h}}{\partial \hat{x}}f(\hat{x}_{t}) + ||\frac{\partial \hat{h}}{\partial \hat{x}}K_{t}c||_{2}\gamma 
 -\frac{1}{2}\mathbf{tr}\left(\sigma^{T}\left(\frac{\partial^{2}\hat{h}}{\partial \hat{x}^{2}}\right)\sigma\right) \\
\label{eq:incomplete-ZCBF-3}
&\leq&\frac{\partial \hat{h}}{\partial \hat{x}}g(\hat{x}_{t})u_{t}
\end{IEEEeqnarray}
where (\ref{eq:incomplete-ZCBF-2}) follows from (\ref{eq:incomplete-ZCBF-1}) and (\ref{eq:incomplete-ZCBF-3}) follows from (\ref{eq:incomplete-ZCBF}). Hence, by Theorem \ref{theorem:ZCBF-safety},  $h(\hat{x}_{t}) > 0$ for all $t$ if $||x_{t}-\hat{x}_{t}||_{2} \leq \gamma$ for all $t$, and thus $Pr\left(h(x_{t}) > 0 \ \forall t\right) \geq (1-\epsilon)$. \qed

\section{CBF-Based Control Policies}
\label{sec:policy}
In what follows, we describe control policies that use stochastic CBFs to provide provable safety guarantees. We consider a case where the goal of the system is to minimize the expected value of a positive-definite quadratic objective function $V(x_{t},u_{t})$. 
 In the complete information case, the controller input $u_{t}$ at time $t$ can be computed as the solution to the quadratic program
\begin{equation}
\label{eq:comp-info-opt}
\begin{array}{ll}
\mbox{minimize}_{u_{t}} & V_{t}(x_{t},u_{t}) \\
\mbox{s.t.} & u_{t} \in \Omega_{t}(x_{t})
\end{array}
\end{equation}
where the set $\Omega_{t}(x_{t})$ is an affine subspace in $u_{t}$. The value of $\Omega_{t}(x_{t})$ depends on whether the RCBF or ZCBF construction is used, as shown in Table \ref{table:linear-constraints}. 

In the incomplete information case, the controller contains  an Extended Kalman Filter, which computes an estimate $\hat{x}_{t}$ of the state $x_{t}$ as a function  of the prior observations $\{y_{\tau} : \tau \in [0,t)\}$. The controller computes each control input $u_{t}$ as a solution to the optimization problem 
\begin{equation}
\label{eq:incomplete-info-opt}
\min{\{V_{t}(\hat{x}_{t},u_{t}) : u_{t} \in \Omega_{t}(\hat{x}_{t})\}}
\end{equation}
where $\Omega_{t}(x_{t})$ is an affine subspace in $u_{t}$. The values of $\Omega_{t}(\hat{x}_{t})$ are shown in Table \ref{table:linear-constraints}.

We observe that these quadratic programs can be extended to describe multiple safety constraints, for example, when the  region $\mathcal{C} = \bigcap_{i=1}^{N}{\{x: h_{i}(x) \geq 0\}}$. This extension can be performed by having a set of linear constraints, one for each safety condition $\{h_{i}(x) \geq 0\}$. There is no guarantee, however, that such a program has a feasible solution $u_{t}$. 

An advantage of the CBF method in the deterministic case is that CBFs can be composed with Control Lyapunov Functions to provide joint guarantees on safety and stability. Such CLFs are defined in the stochastic setting as follows.
\begin{proposition}[\cite{florchinger1997feedback}]
\label{prop:CLF}
Suppose there exists a function $V: \mathbb{R}^{n} \rightarrow \mathbb{R}$ such that, for every $x$, there exists $u$ satisfying 
\begin{equation}
\label{eq:CLF}
\frac{\partial V}{\partial x}(f(x) + g(x)u) + \mathbf{tr}\left(\sigma^{T}\frac{\partial^{2}V}{\partial x^{2}}\sigma\right) \leq 0
\end{equation}
If $u_{t}$ is chosen to satisfy (\ref{eq:CLF}) at each time $t$, then $0$ is stochastically asymptotically stable.
\end{proposition}
Proposition \ref{prop:CLF} implies that stability requirements can be incorporated as a linear constraint on the optimization-based control. Hence, if the control input can be chosen at each time $t$ to jointly satisfy the appropriate CBF constraint of Table \ref{table:linear-constraints} and the CLF constraint (\ref{eq:CLF}), then the system is guaranteed to asymptotically approach the desired operating point while remaining safe for all time. 

\begin{figure*}
\centering
$\begin{array}{ccc}
\includegraphics[width=2in]{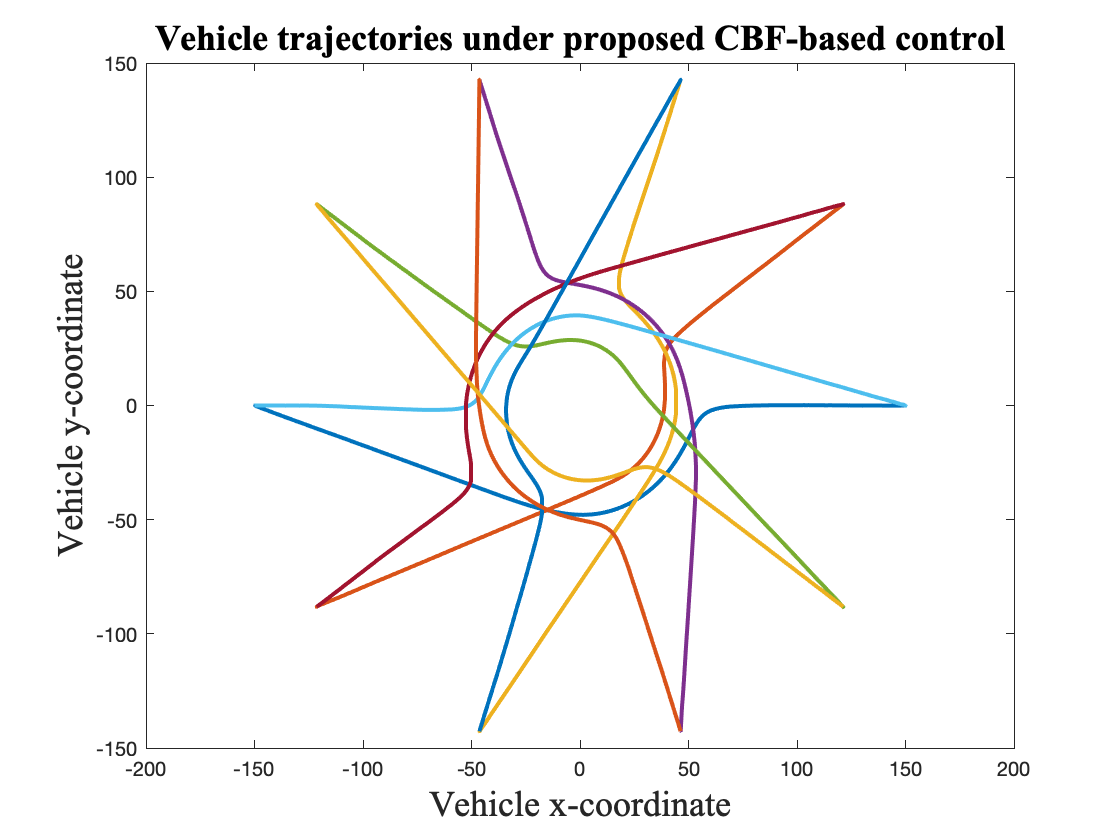} &
\includegraphics[width=2in]{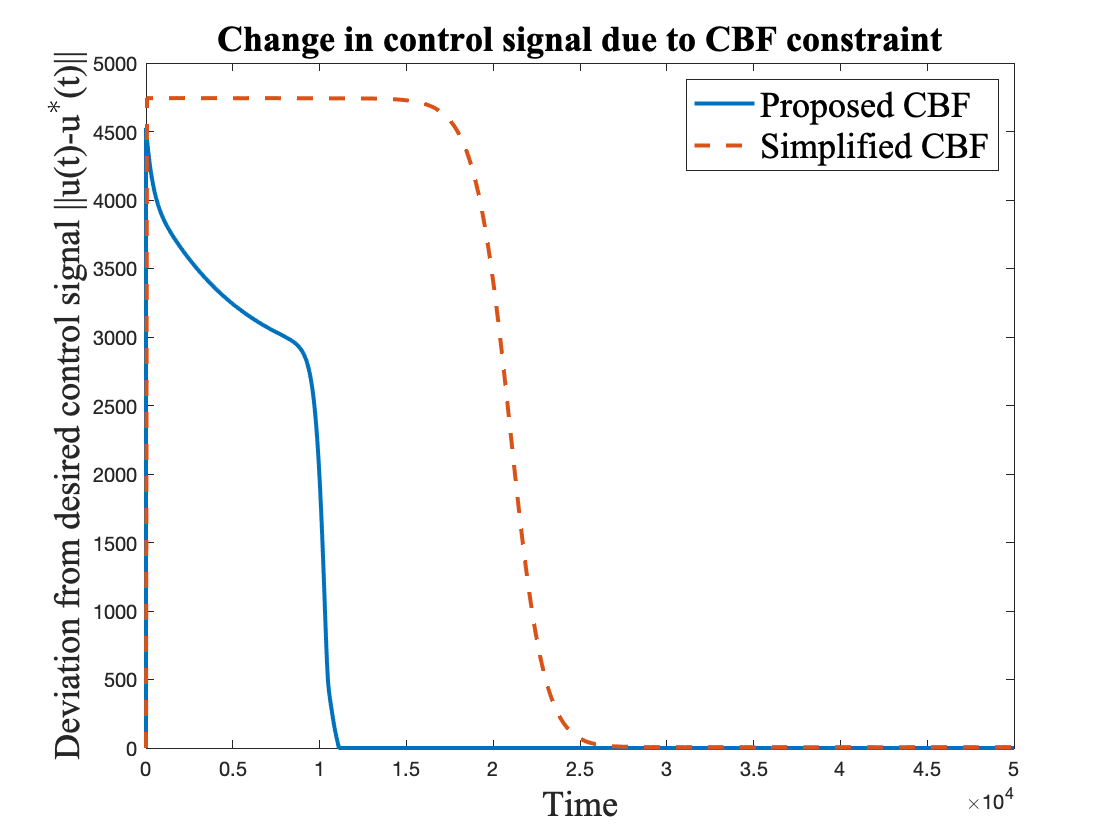} & 
\includegraphics[width=2in]{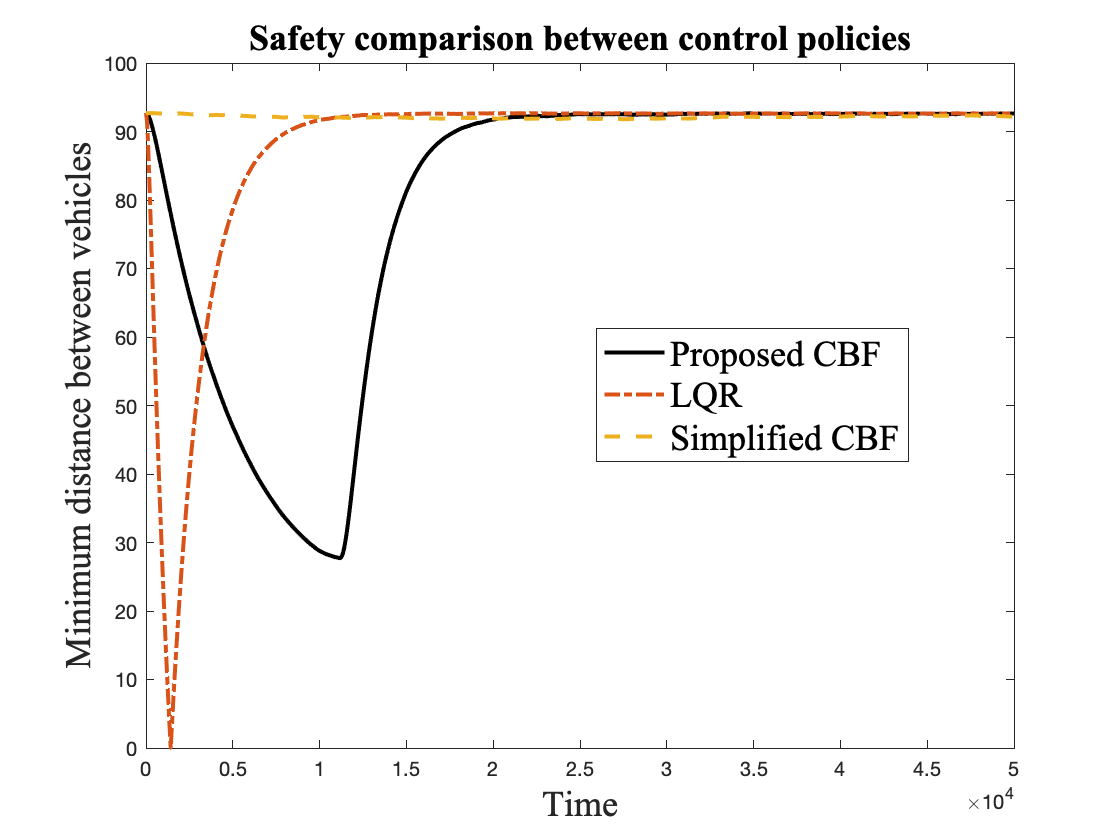} 
\\
\mbox{(a)} & \mbox{(b)} & \mbox{(c)}
\end{array}$
\caption{Evaluation of stochastic CBFs using multi-agent collision avoidance. Ten agents are initially placed equidistant on a circle and must reach the opposite point without colliding in the presence of process and measurement noise. We compare our stochastic ZCBF approach with a proportional linear control law that does not incorporate safety, as well as a simplified CBF that uses the deterministic CBF of \cite{ames2016control} on the estimated state. (a) The agent trajectories. Each agent has double-integrator dynamics and avoids the center of the circle to prevent collisions. (b) Comparison of the deviation between the control action chosen and the desired action under a linear control law. The simplified CBF leads to higher deviation compared to our proposed stochastic ZCBF. (c) The minimum distance between vehicles under the linear control law and two CBF-based approaches. The linear controller leads to safety violations, while the CBF-based approaches avoid collisions.}
\label{fig:sim}  
\end{figure*}

\section{Numerical Study}
\label{sec:simulation}
We performed a numerical study of a multi-agent collision avoidance scenario using Matlab. Our case study is based on \cite{borrmann2015control}. We considered a set of $n$ agents, indexed $i=1,\ldots,n$, where agent $i$ has position and velocity $[p_{t}]_{i}$ and $[v_{t}]_{i}$ with $d[p_{t}]_{i} = [v_{t}]_{i} \ dt + \sigma_{p} \ dW_{t}$ and $d[v_{t}]_{i}= [u_{t}]_{i} \ dt + \sigma_{v} \ dW_{t}$. The agents were uniformly placed on a circle of radius $\rho = 150$, with each agent attempting to travel to the opposite point on the circle while avoiding collisions. Each pair of agents $(i,j)$ had a safety constraint $$h_{ij} = ||p_{i}-p_{j}||_{2} - D_{s} \geq 0,$$ where $D_{s} = 10$. The sensor measurements satisfied $dy_{t} = x_{t} \ dt + \nu \ dW_{t}$ where $\nu = I$. We set $\sigma_{p} = \sigma_{v}= I$. The cost function to be minimized was equal to $||u - \bar{u}||_{2}$, where $\bar{u}$ is obtained using a linear control law $\bar{u}_{t} = -k_{1}([p_{t}]_{i}-[r]_{i}) - k_{2}[v_{t}]_{i}$. The CBF constraints were constructed using the method for high relative degree systems introduced in Section \ref{subsec:high-degree-complete}.

The agent trajectories are shown in Figure \ref{fig:sim}(a). Each agent moves to reach the desired destination while avoiding collisions. We observe that all agents avoid traversing the center in order to minimize collisions. Fig. \ref{fig:sim}(b) shows the gap between the desired and actual control input over time. The proposed stochastic ZCBF led to a reduced deviation from the desired control input compared to the simplified CBF.

The minimum distances achieved by the three policies are shown in Fig. \ref{fig:sim}(c). The linear control law leads to safety violations as the agents attempt to reach their desired final positions while disregarding safety. The CBF-based approaches both avoid safety violations, with the stochastic ZCBF approaching the unsafe region before recovering to maintain a safe distance and still converging to the desired final position.
\section{Conclusion}
\label{sec:conclusion}
This paper developed a framework for safe control of stochastic systems via Control Barrier Functions. We considered two scenarios, namely, complete information in which the true  state is known to the controller at each time, and incomplete information in which the controller only has access to sensor measurements that are corrupted by Gaussian noise. For each case, we constructed Reciprocal and Zero CBFs, which ensure that the system remains safe provided that the CBF is finite (RCBF) or nonnegative (ZCBF).  We proved that both constructions guarantee safety with probability 1 in the complete information case, and provide stochastic safety guarantees that depend on the estimation accuracy in the incomplete information case. We proposed control policies that ensure safety and stability by solving quadratic programs containing CBFs and stochastic Control Lyapunov Functions (CLFs) at each time step.  We evaluated our approach through a numerical simulation on a multi-agent collision avoidance scenario. Future work will consider techniques for more general high relative-degree systems, as well as systems that are not affine in the control input. Another direction for future work consists of exploring the distinctions between RCBF- and ZCBF-based control policies. For example, in the deterministic case, the fact that the ZCBF is well-defined even outside the safe region can be used to design controllers that converge  to the safe region if the initial state is outside the safe region. Generalizing such results to the stochastic setting remains an open problem.


\bibliographystyle{plain}        
\bibliography{autosam}           



\end{document}